\documentclass[11pt]{article}

\usepackage[T1]{fontenc}
\usepackage[utf8]{inputenc}
\usepackage{lmodern}
\usepackage{amsmath,amssymb,amsthm,mathtools}
\usepackage[margin=1in]{geometry}
\usepackage{microtype}
\usepackage{xcolor}
\usepackage[
  colorlinks=true,
  linkcolor=blue,
  citecolor=red,
  urlcolor=blue,
  linktoc=all
]{hyperref}

\hypersetup{
  pdftitle={The curse of dimensionality for geometric L1-discrepancies with nonnegative weights},
  pdfauthor={Josef Dick and Esra Sengelen Sevim},
  pdfsubject={Geometric L1-discrepancy and the curse of dimensionality},
  pdfkeywords={L1-discrepancy, curse of dimensionality, nonnegative weights, information complexity}
}

\numberwithin{equation}{section}

\theoremstyle{plain}
\newtheorem{theorem}{Theorem}[section]

\newtheorem{lemma}[theorem]{Lemma}

\theoremstyle{definition}

\newtheorem{remark}[theorem]{Remark}

\newcommand{\R}{\mathbb R}
\newcommand{\N}{\mathbb N}
\newcommand{\T}{\mathbb T}
\newcommand{\1}{\mathbf 1}
\newcommand{\dd}{\,\mathrm d}
\newcommand{\disc}{\operatorname{disc}}
\newcommand{\Estar}{E^{\mathrm{star}}}
\newcommand{\Eext}{E^{\mathrm{ext}}}
\newcommand{\Eper}{E^{\mathrm{per}}}

\newcommand{\doi}[1]{\href{https://doi.org/#1}{\nolinkurl{doi:#1}}}

\title{The curse of dimensionality for geometric
\texorpdfstring{$L_1$}{L1}-discrepancies with nonnegative weights}
\author{
  Josef Dick\thanks{josef.dick@unsw.edu.au}\\[2pt]
  School of Mathematics and Statistics,\\
  UNSW Sydney, Sydney NSW 2052, Australia
  \and
  Esra Sengelen Sevim\thanks{esra.sevim@nisantasi.edu.tr}\\[2pt]
  Istanbul Ni\c{s}anta\c{s}\i{} University,\\ Istanbul, Turkey}
\date{}

\begin{document}

\maketitle

\begin{abstract}
We prove that the normalized star, extreme, and periodic
$L_1$-discrepancies with arbitrary nonnegative weights suffer from the
curse of dimensionality.  A single elementary lemma is used in all three cases:
after normalizing the target and every one-dimensional kernel to probability
densities, a uniform square-root-affinity gap tensorizes over the coordinates.
For every $d\in\N$ and $0<\varepsilon<1$, the normalized information complexities have lower bounds of the
form
\[
 \frac{(1-\varepsilon)^2}{1+\varepsilon}\,q^d,
\]
with
\[
 q_{\mathrm{star}}=\frac12+\frac1{\sqrt3},\qquad
 q_{\mathrm{ext}}=\frac{75(9+4\sqrt2)}{784},\qquad
 q_{\mathrm{per}}=\frac{25}{16}.
\]
All three constants are larger than one.  The argument uses only
nonnegativity, product integration, the elementary inequalities for the
square root, and Cauchy--Schwarz.
\end{abstract}



\section{Introduction}

Geometric discrepancy measures how well a finite weighted set reproduces
the volumes of a prescribed family of test boxes. It originates in the
study of uniform distribution and is a central tool in quasi-Monte Carlo
(QMC) integration \cite{DickPillichshammer2010}. Its inverse asks how many
nodes are needed to reduce the discrepancy by a fixed factor relative to
the empty rule. Normalization by this initial error is standard in
information-based complexity; see the three-volume monograph
\cite{NovakWozniakowski2008,NovakWozniakowski2010,NovakWozniakowski2012}.
The relationship between discrepancy, integration, and tractability is
surveyed in \cite{Hinrichs2013}. In the star-discrepancy setting,
nonnegative node weights correspond to positive quadrature formulas in
suitable anchored Sobolev spaces
\cite{NovakPillichshammer2025PAMS,NovakPillichshammer2026Survey}.

There are two different asymptotic regimes in discrepancy theory. For
fixed dimension and increasing numbers of points, Roth's foundational
$L_2$ lower bound \cite{Roth1954} initiated the modern theory of
$L_p$-discrepancy. Average and probabilistic versions describe the behavior
of random point sets. Relevant results include the estimates for average
extreme discrepancy in \cite{Gnewuch2005}, the asymptotic formulas of
Steinerberger \cite{Steinerberger2010} and Hinrichs--Weyhausen
\cite{HinrichsWeyhausen2012}, and the analysis of generalized random point
sets and importance sampling in \cite{NovakPillichshammer2026Random}.
Such results provide useful context, but average-case estimates alone do
not determine the tractability of the optimal deterministic discrepancy.
Our concern here is the dependence of the normalized information
complexity on the dimension.

For the classical star discrepancy, that is, the $L_\infty$-norm of the
anchored discrepancy function, Heinrich, Novak, Wasilkowski, and
Wo\'zniakowski proved that the equal-weight information complexity satisfies
\[
 N^{\mathrm{star},\mathrm{qmc}}_\infty(\varepsilon,d)
 \le C d\varepsilon^{-2}
\]
with an absolute constant $C>0$ \cite{HeinrichEtAl2001}.
Hinrichs proved a corresponding linear lower bound in $d$ for fixed
sufficiently small $\varepsilon$ \cite{Hinrichs2004}; see also the
elementary proof in \cite{Steinerberger2023}. Aistleitner and Hofer
obtained explicit probabilistic bounds for Monte Carlo point sets
\cite{AistleitnerHofer2014}, while Doerr proved that, for $d\le N$,
independent uniformly distributed points have expected star discrepancy
of order $\sqrt{d/N}$ \cite{Doerr2014}. The extreme and periodic
$L_\infty$-discrepancies are also polynomially tractable; see
\cite{Gnewuch2005,NovakPillichshammer2026Survey}.

The situation is different for finite $p>1$. Early intractability results
for integration and discrepancy were obtained by Wo\'zniakowski and by
Novak--Wo\'zniakowski
\cite{Wozniakowski2000,NovakWozniakowski2001}. Following the results in
\cite{NovakPillichshammer2023}, Novak and Pillichshammer established the
curse of dimensionality for the normalized star $L_p$-discrepancy with
nonnegative weights for every $p\in(1,\infty)$
\cite{NovakPillichshammer2025PAMS}. Their general tensor-product
lower-bound method is developed in \cite{NovakPillichshammer2025Tensor}.
The corresponding results for the extreme and periodic geometries appear
in \cite{NovakPillichshammer2026Extreme,NovakPillichshammer2026Survey}.
The weighted $L_p$-discrepancy, in which coordinate weights model unequal
importance of variables, has a separate tractability theory
\cite{SloanWozniakowski1998,NovakPillichshammer2025Weighted}. These
coordinate weights should be distinguished from the nonnegative node
weights considered here. Other norms can lead to different behavior:
discrepancy in certain Orlicz norms is polynomially or weakly tractable
\cite{DickHinrichsPillichshammerProchno2020}, whereas BMO discrepancy
suffers from the curse of dimensionality \cite{Pillichshammer2023}.

The endpoint $p=1$ was not covered by these finite-$p$ lower-bound
arguments. The analysis of random points in
\cite{NovakPillichshammer2026Random} retains exponential dimension
dependence even after optimizing the product sampling density and its
associated nonnegative importance weights. This provides motivation for
a deterministic lower bound, but does not itself prove one. The survey
\cite{NovakPillichshammer2026Survey} explicitly posed the $L_1$ cases for
the star, extreme, and periodic discrepancies as open problems; the
extreme case was also left open in
\cite{NovakPillichshammer2026Extreme}.

In this paper we resolve all three endpoints for arbitrary nonnegative
node weights. A single elementary product-kernel argument gives explicit
exponential lower bounds for the normalized information complexities.
The proof uses a volume-biased probability measure and a square-root
moment. After normalization, the one-dimensional quantity is the
Hellinger affinity between the target density and a kernel density. A
uniform upper bound strictly smaller than one for this affinity
multiplies over the coordinates and yields the desired lower bound.

The proof has two parts. Section~\ref{sec:product} proves the abstract
lower bound for a nonnegative sum of product kernels. Relative $L_1$
accuracy controls the total normalized mass of the sum, while the
one-dimensional affinity bound tensorizes over the coordinates. In
Section~\ref{sec:affinities} we compute the required affinity for the
star, extreme, and periodic test classes. Section~\ref{sec:main-proof}
passes from bounds for individual rules to minimal discrepancies and
information complexities.

The method applies to arbitrary nonnegative coefficients; their sum is
not assumed to be one. Consequently, the same lower bounds hold for
classical equal-weight point sets. Nonnegativity is essential: it is used
when taking the square root of the finite sum, and the argument does not
apply to signed weights.

\section{Discrepancies and normalized information complexity}
\label{sec:definitions}

Throughout, $\N=\{1,2,\ldots\}$, $\N_0=\{0,1,2,\ldots\}$,
$d\in\N$. In the discrepancy definitions below, integration is with
respect to Lebesgue measure. A weighted rule with at most $N$ terms consists of nodes
\(x^{(1)},\ldots,x^{(N)}\in[0,1)^d\)
and coefficients $a_1,\ldots,a_N\ge0$.  Repeated nodes and zero coefficients
are allowed.  This only enlarges the usual class and therefore makes a lower
bound stronger.  The equal-weight, or QMC, class is obtained by requiring
 $a_i=1/N$.

\subsection{Star discrepancy}

For $t=(t_1,\ldots,t_d)\in[0,1]^d$, put
\([0,t)=\prod_{j=1}^d[0,t_j)\).  Define
\begin{equation}\label{eq:star-error}
 \Estar_{N,d}(P,a)
 =\int_{[0,1]^d}
 \left|
   \sum_{i=1}^N a_i\1_{[0,t)}(x^{(i)})
   -\prod_{j=1}^d t_j
 \right|\dd t.
\end{equation}
The empty rule has error
\begin{equation}\label{eq:star-initial}
 \Estar_{0,d}=\left(\int_0^1t\dd t\right)^d=2^{-d}.
\end{equation}

\subsection{Extreme discrepancy}

Let
\[
 \mathcal D=\{(u,v)\in[0,1]^2: u\le v\},
 \qquad
 \mathcal D^d=\mathcal D\times\cdots\times\mathcal D,
\]
and, for $(u,v)\in\mathcal D^d$, let
\([u,v)=\prod_{j=1}^d[u_j,v_j)\).  Define
\begin{equation}\label{eq:ext-error}
 \Eext_{N,d}(P,a)
 =\int_{\mathcal D^d}
 \left|
   \sum_{i=1}^N a_i\1_{[u,v)}(x^{(i)})
   -\prod_{j=1}^d(v_j-u_j)
 \right|\dd(u,v).
\end{equation}
Its initial error is
\begin{equation}\label{eq:ext-initial}
 \Eext_{0,d}
 =\left(\int_0^1\int_u^1(v-u)\dd v\dd u\right)^d
 =6^{-d}.
\end{equation}

\subsection{Periodic discrepancy}

Identify $\T=[0,1)$ with the circle $\R/\mathbb Z$.  For
$u\in\T$ and $0\le\ell<1$, let
\[
 I(u,\ell)=\{u+s\pmod1:0\le s<\ell\}.
\]
For $u\in\T^d$ and $\ell\in[0,1)^d$, put
\(B(u,\ell)=\prod_{j=1}^d I(u_j,\ell_j)\), and define
\begin{equation}\label{eq:per-error}
 \Eper_{N,d}(P,a)
 =\int_{\T^d}\int_{[0,1)^d}
 \left|
   \sum_{i=1}^N a_i\1_{B(u,\ell)}(x^{(i)})
   -\prod_{j=1}^d\ell_j
 \right|\dd\ell\dd u.
\end{equation}
This length parametrization is equivalent to the usual endpoint
parametrization: the change of variables
\((u,v)\mapsto(u,\ell)=(u,v-u\pmod1)\)
preserves Lebesgue measure on $\T^2$.  Indeed, for each fixed $u$, the map
$v\mapsto v-u\pmod1$ is a translation of the circle and hence preserves
one-dimensional Lebesgue measure; Fubini's theorem gives the assertion.
The empty-rule error is
\begin{equation}\label{eq:per-initial}
 \Eper_{0,d}
 =\left(\int_0^1\int_0^1\ell\dd\ell\dd u\right)^d
 =2^{-d}.
\end{equation}

\subsection{Minimal errors and the curse of dimensionality}

For
\(\mathcal C\in\{\mathrm{star},\mathrm{ext},\mathrm{per}\}\),
let
\[
 \disc^{\mathcal C,+}_1(N,d)
 =\inf_{P,a} E^{\mathcal C}_{N,d}(P,a),
\]
where the infimum is over all rules with $N$ nonnegative coefficients.  We
use the initial errors in \eqref{eq:star-initial}, \eqref{eq:ext-initial},
and \eqref{eq:per-initial} when $N=0$.  For $0<\varepsilon<1$, define
\begin{equation}\label{eq:complexity}
 N^{\mathcal C,+}_1(\varepsilon,d)
 =\min\left\{N\in\N_0:
 \disc^{\mathcal C,+}_1(N,d)
 \le\varepsilon E^{\mathcal C}_{0,d}\right\},
\end{equation}
with $\min\varnothing=\infty$.  The QMC quantities
\(\disc^{\mathcal C,\mathrm{qmc}}_1\)
and
\(N^{\mathcal C,\mathrm{qmc}}_1\)
are defined by restricting to equal weights.  Since this is a subclass,
\begin{equation}\label{eq:qmc-subclass}
 N^{\mathcal C,\mathrm{qmc}}_1(\varepsilon,d)
 \ge N^{\mathcal C,+}_1(\varepsilon,d).
\end{equation}

We say that a problem suffers from the curse of dimensionality if, for some
fixed $0<\varepsilon<1$, there is a constant $c>1$ such that its information
complexity is at least $c^d$ for infinitely many $d$.  The estimates below
give this conclusion for every fixed $\varepsilon\in(0,1)$ and all
sufficiently large $d$.

The main result is the following.

\begin{theorem}[Main theorem]\label{thm:main}
For every $d\in\N$ and $0<\varepsilon<1$,
\begin{align}
 N^{\mathrm{star},+}_1(\varepsilon,d)
 &\ge \frac{(1-\varepsilon)^2}{1+\varepsilon}
 \left(\frac12+\frac1{\sqrt3}\right)^d,
 \label{eq:main-star}\\
 N^{\mathrm{ext},+}_1(\varepsilon,d)
 &\ge \frac{(1-\varepsilon)^2}{1+\varepsilon}
 \left(\frac{75(9+4\sqrt2)}{784}\right)^d, \label{eq:main-ext}
 \\
 N^{\mathrm{per},+}_1(\varepsilon,d)
 &\ge \frac{(1-\varepsilon)^2}{1+\varepsilon}
 \left(\frac{25}{16}\right)^d.
 \label{eq:main-per}
\end{align}
The same three bounds hold with the superscript $+$ replaced by
$\mathrm{qmc}$.  Consequently, the normalized star, extreme, and periodic
$L_1$-discrepancies with nonnegative weights, and their equal-weight
subclasses, suffer from the curse of dimensionality.
\end{theorem}

\section{An elementary product-kernel lower bound}
\label{sec:product}

We isolate the only abstract argument used in the paper.  The statement is
given for general measure spaces so that the three geometric applications
become direct substitutions.

Let $(\Omega,\Sigma,\mu)$ be a $\sigma$-finite measure space.  Let
$h:\Omega\to[0,\infty)$ be measurable, positive $\mu$-almost everywhere,
and satisfy
\begin{equation*}
 0<I:=\int_\Omega h\dd\mu<\infty.
\end{equation*}
Let $\{k_x:x\in X\}$ be nonnegative measurable kernels with finite masses
\[
 m(x)=\int_\Omega k_x\dd\mu<\infty.
\]
For $m(x)>0$, define the probability measure $\nu$ and the normalized kernel
density $f_x$ by
\begin{equation}\label{eq:normalized-density}
 \dd\nu=\frac hI\dd\mu,
 \qquad
 f_x=\frac{I k_x}{m(x)h}.
\end{equation}
The value of $f_x$ on the null set $\{h=0\}$ is irrelevant.  Directly from
\eqref{eq:normalized-density},
\begin{equation}\label{eq:density-and-affinity}
 \int_\Omega f_x\dd\nu=1,
 \qquad
 A(x):=\int_\Omega\sqrt{f_x}\dd\nu
 =\frac1{\sqrt{I m(x)}}
   \int_\Omega\sqrt{h k_x}\dd\mu.
\end{equation}

For $d\in\N$ and $\omega=(\omega_1,\ldots,\omega_d)\in\Omega^d$, set
\[
 H_d(\omega)=\prod_{j=1}^dh(\omega_j),
 \qquad
 K_y(\omega)=\prod_{j=1}^dk_{y_j}(\omega_j),
 \quad y\in X^d.
\]

\begin{lemma}[Square-root product-kernel bound]
\label{lem:product}
Assume that there is a number $0<\rho<1$ such that
\begin{equation}\label{eq:rho-assumption}
 A(x)\le\rho
 \qquad\text{whenever }m(x)>0.
\end{equation}
Let
\[
 G=\sum_{i=1}^M a_iK_{y^{(i)}},
 \qquad
 a_i>0, \quad m(y^{(i)}_j)>0,
\]
and suppose that
\begin{equation}\label{eq:relative-error-assumption}
 \|H_d-G\|_{L_1(\mu^d)}\le\varepsilon I^d,
 \qquad 0<\varepsilon<1.
\end{equation}
Then
\begin{equation}\label{eq:abstract-lower-bound}
 M\ge\frac{(1-\varepsilon)^2}{1+\varepsilon}\,\rho^{-2d}.
\end{equation}
\end{lemma}

\begin{proof}
Let $\nu_d=\nu^{d}$.  On the set where $H_d>0$, define
\(Z=G/H_d\).  This set has full $\nu_d$-measure.  For $i=1,\ldots,M$, put
\[
 b_i=a_i\prod_{j=1}^d\frac{m(y_j^{(i)})}{I},
 \qquad
 F_i(\omega)=\prod_{j=1}^df_{y_j^{(i)}}(\omega_j).
\]
Then
\begin{equation}\label{eq:Z-sum}
 Z=\sum_{i=1}^M b_iF_i
 \quad\nu_d\text{-almost surely},
 \qquad
 \int_{\Omega^d}F_i\dd\nu_d=1.
\end{equation}
Moreover, product integration and \eqref{eq:rho-assumption} give
\begin{equation}\label{eq:F-affinity}
 \int_{\Omega^d}\sqrt{F_i}\dd\nu_d
 =\prod_{j=1}^d A(y_j^{(i)})
 \le\rho^d.
\end{equation}

The density of $\nu_d$ with respect to $\mu^d$ is $H_d/I^d$.  Hence
\begin{equation}\label{eq:normalized-L1}
 \int_{\Omega^d}|1-Z|\dd\nu_d
 =\frac{\|H_d-G\|_{L_1(\mu^d)}}{I^d}
 \le\varepsilon.
\end{equation}
Let $B=\sum_{i=1}^M b_i$.  By \eqref{eq:Z-sum},
\(B=\int Z\dd\nu_d\), and therefore
\begin{equation}\label{eq:B-upper}
 B\le1+\int|Z-1|\dd\nu_d\le1+\varepsilon.
\end{equation}

For every $z\ge0$,
\begin{equation}\label{eq:sqrt-distance}
 \sqrt z\ge1-|z-1|.
\end{equation}
Indeed,
\(|1-\sqrt z|\le|1-z|\)
because
\(|1-z|=|1-\sqrt z|(1+\sqrt z)\).
Combining \eqref{eq:sqrt-distance} and \eqref{eq:normalized-L1} yields
\begin{equation}\label{eq:sqrt-lower}
 \int\sqrt Z\dd\nu_d\ge1-\varepsilon.
\end{equation}

For nonnegative numbers,
\(\sqrt{r+s}\le\sqrt r+\sqrt s\).  Applying this inequality repeatedly to
\eqref{eq:Z-sum}, then using \eqref{eq:F-affinity}, Cauchy--Schwarz, and
\eqref{eq:B-upper}, we obtain
\begin{align}
 \int\sqrt Z\dd\nu_d
 &\le\sum_{i=1}^M\sqrt{b_i}
       \int\sqrt{F_i}\dd\nu_d\notag\\
 &\le\rho^d\sum_{i=1}^M\sqrt{b_i}
 \le\rho^d\sqrt{M\sum_{i=1}^Mb_i}
 \le\rho^d\sqrt{M(1+\varepsilon)}.
 \label{eq:sqrt-upper}
\end{align}
Together, \eqref{eq:sqrt-lower} and \eqref{eq:sqrt-upper} imply
\[
 1-\varepsilon\le\rho^d\sqrt{M(1+\varepsilon)}.
\]
Squaring and rearranging proves \eqref{eq:abstract-lower-bound}.
\end{proof}

\begin{remark}\label{rem:zero-terms}
The positivity assumptions on the displayed terms cause no loss of
generality.  A term with coefficient zero may be deleted.  If $m(x)=0$,
then the nonnegative kernel $k_x$ vanishes almost everywhere.  Hence a
product kernel with a zero-mass coordinate also vanishes almost everywhere
and may be deleted.  If the original rule contains at most $N$ terms, the
remaining number $M$ satisfies $M\le N$.  Since $\varepsilon<1$, at least one
term remains whenever \eqref{eq:relative-error-assumption} holds.
\end{remark}

\section{The three one-dimensional affinity bounds}
\label{sec:affinities}

We now calculate $I$, $m$, and the maximum of the affinity $A$ in
\eqref{eq:density-and-affinity}.  These are the only geometry-specific
steps.

\subsection{Star boxes}
\label{subsec:star-affinity}

Take
\[
 \Omega=[0,1],\qquad h(t)=t,\qquad
 k_x(t)=\1_{\{x<t\}},\quad 0\le x<1.
\]
Then
\begin{equation*}
 I=\frac12, \qquad m(x)=1-x.
\end{equation*}
By \eqref{eq:density-and-affinity},
\begin{equation}\label{eq:star-affinity}
 A_{\mathrm{star}}(x)
 =\frac{2\sqrt2}{3\sqrt{1-x}}(1-x^{3/2}).
\end{equation}

\begin{lemma}\label{lem:star-affinity}
For $0\le x<1$,
\begin{equation*}
 A_{\mathrm{star}}(x)
 \le\rho_{\mathrm{star}}:=\sqrt{4\sqrt3-6}<1.
\end{equation*}
Equality holds at $x=1-\sqrt3/2$.  Consequently,
\begin{equation}\label{eq:star-base}
 \rho_{\mathrm{star}}^{-2}
 =\frac1{4\sqrt3-6}
 =\frac12+\frac1{\sqrt3}.
\end{equation}
\end{lemma}

\begin{proof}
Put $u=\sqrt x$.  Squaring \eqref{eq:star-affinity} and factoring
$1-u^3$ and $1-u^2$ gives
\begin{equation}\label{eq:star-H}
 A_{\mathrm{star}}(u^2)^2
 =\frac89 H(u),
 \qquad
 H(u)=\frac{(1-u)(1+u+u^2)^2}{1+u},
 \quad 0\le u\le1.
\end{equation}
The value at $u=1$ is understood by continuity.  Direct differentiation
gives
\[
 H'(u)
 =-\frac{2u(1+u+u^2)(2u^2+2u-1)}{(1+u)^2}.
\]
The only sign change in $(0,1)$ occurs at
\[
 u_0=\frac{\sqrt3-1}{2}.
\]
Thus $H$ increases on $(0,u_0)$ and decreases on $(u_0,1)$.  Since
\[
 1+u_0+u_0^2=\frac32,
 \qquad
 \frac{1-u_0}{1+u_0}=2\sqrt3-3,
\]
we have
\[
 H(u_0)=\frac94(2\sqrt3-3).
\]
Equation \eqref{eq:star-H} now gives
\(A_{\mathrm{star}}(u_0^2)^2=4\sqrt3-6\).
Finally,
\(u_0^2=1-\sqrt3/2\), and rationalizing the denominator proves
\eqref{eq:star-base}.
\end{proof}

\subsection{Extreme boxes}
\label{subsec:extreme-affinity}

Take
\[
 \Omega=\mathcal D, \qquad h(u,v)=v-u, \qquad
 k_z(u,v)=\1_{\{u\le z<v\}},\quad 0\le z\le1.
\]
Direct integration gives
\begin{equation}\label{eq:ext-Im}
 I=\frac16, \qquad m(z)=z(1-z).
\end{equation}
The endpoint kernels have mass zero.  For $0<z<1$, define
\begin{equation}\label{eq:M-def}
 M(z)=\int_0^1\int_0^1
 \sqrt{zx+(1-z)y}\dd x\dd y.
\end{equation}

\begin{lemma}\label{lem:ext-affinity}
For $0<z<1$,
\begin{equation}\label{eq:ext-affinity}
 A_{\mathrm{ext}}(z)
 =\sqrt{6z(1-z)}\,M(z)
 \le\rho_{\mathrm{ext}}
 :=\frac4{15}\sqrt{\frac32}(4-\sqrt2)<1.
\end{equation}
Equality holds at $z=1/2$.  Consequently,
\begin{equation}\label{eq:ext-base}
 \rho_{\mathrm{ext}}^{-2}
 =\frac{75}{8(4-\sqrt2)^2}
 =\frac{75(9+4\sqrt2)}{784}.
\end{equation}
\end{lemma}

\begin{proof}
By \eqref{eq:density-and-affinity} and \eqref{eq:ext-Im},
\begin{align*}
 A_{\mathrm{ext}}(z)
 &=\sqrt{\frac6{z(1-z)}}
   \int_0^z\int_z^1\sqrt{v-u}\dd v\dd u.
\end{align*}
Use the substitutions
\[
 u=z(1-x),\qquad v=z+(1-z)y.
\]
The Jacobian is $z(1-z)$ and
\(v-u=zx+(1-z)y\).  This proves the identity in
\eqref{eq:ext-affinity}.

For fixed $x,y\in[0,1]$, the function
\(z\mapsto\sqrt{zx+(1-z)y}\)
is concave.  Therefore $M$ is concave.  Interchanging $x$ and $y$ in
\eqref{eq:M-def} shows that $M(z)=M(1-z)$.  Hence
\[
 M(1/2)
 \ge\frac{M(z)+M(1-z)}2=M(z).
\]
Since also $z(1-z)\le1/4$, both nonnegative factors in
\eqref{eq:ext-affinity} are maximized at $z=1/2$.

It remains to evaluate one elementary double integral.  Integrating first
in $y$ and then in $x$ gives
\begin{align*}
 M(1/2)
 &=2^{-1/2}\int_0^1\int_0^1\sqrt{x+y}\dd y\dd x\\
 &=\frac{4(4-\sqrt2)}{15}.
\end{align*}
Thus
\[
 A_{\mathrm{ext}}(1/2)
 =\sqrt{\frac32}\,\frac{4(4-\sqrt2)}{15},
\]
which is the stated value of $\rho_{\mathrm{ext}}$.  Squaring gives
\[
 \rho_{\mathrm{ext}}^2
 =\frac{8(4-\sqrt2)^2}{75}<1.
\]
Taking the reciprocal and rationalizing proves \eqref{eq:ext-base}.
\end{proof}

\subsection{Periodic boxes}
\label{subsec:periodic-affinity}

It is convenient to use the starting point and length as the two
one-dimensional parameters.  Take
\[
 \Omega=\T\times[0,1),\qquad
 h(u,\ell)=\ell, \qquad
 k_z(u,\ell)=\1_{I(u,\ell)}(z),\quad z\in\T.
\]
For a fixed $z$ and $\ell$, the set of starting points $u$ for which
$z\in I(u,\ell)$ is a circle interval of length $\ell$.  Hence
\begin{equation}\label{eq:per-Im}
 I=\frac12, \qquad
 m(z)=\int_0^1\ell\dd\ell=\frac12.
\end{equation}

\begin{lemma}\label{lem:per-affinity}
For every $z\in\T$,
\begin{equation*}
 A_{\mathrm{per}}(z)=\frac45=:\rho_{\mathrm{per}}.
\end{equation*}
Consequently,
\begin{equation}\label{eq:per-base}
 \rho_{\mathrm{per}}^{-2}=\frac{25}{16}.
\end{equation}
\end{lemma}

\begin{proof}
By \eqref{eq:density-and-affinity} and \eqref{eq:per-Im},
\begin{align*}
 A_{\mathrm{per}}(z)
 &=2\int_0^1\int_\T
   \sqrt\ell\,\1_{I(u,\ell)}(z)\dd u\dd\ell\\
 &=2\int_0^1\ell^{3/2}\dd\ell
 =\frac45.
\end{align*}
The second equality uses that the relevant set of $u$ has measure $\ell$.
Equation \eqref{eq:per-base} follows immediately.
\end{proof}

\section{Proof of the main theorem}
\label{sec:main-proof}

We first prove a bound for an individual nonnegative rule.  In each of the
three cases, the discrepancy integrand has the product-kernel form from
Section~\ref{sec:product}.

For the star discrepancy, take $h$ and $k_x$ from
Section~\ref{subsec:star-affinity}.  Up to parameter sets of measure zero,
\[
 \prod_{j=1}^dk_{x_j^{(i)}}(t_j)
 =\1_{[0,t)}(x^{(i)}),
 \qquad
 \prod_{j=1}^dh(t_j)=\prod_{j=1}^dt_j.
\]
Thus \eqref{eq:star-error} is exactly
\(\|H_d-G\|_{L_1}\), and $I^d=2^{-d}$ is the initial error.

For the extreme discrepancy, use $h$ and $k_z$ from
Section~\ref{subsec:extreme-affinity}.  Then
\[
 \prod_{j=1}^dk_{x_j^{(i)}}(u_j,v_j)
 =\1_{[u,v)}(x^{(i)}),
 \qquad
 \prod_{j=1}^dh(u_j,v_j)=\prod_{j=1}^d(v_j-u_j).
\]
Hence \eqref{eq:ext-error} is $\|H_d-G\|_{L_1}$ and
$I^d=6^{-d}$ is the initial error.

For the periodic discrepancy, use $h$ and $k_z$ from
Section~\ref{subsec:periodic-affinity}.  Then
\[
 \prod_{j=1}^dk_{x_j^{(i)}}(u_j,\ell_j)
 =\1_{B(u,\ell)}(x^{(i)}),
 \qquad
 \prod_{j=1}^dh(u_j,\ell_j)=\prod_{j=1}^d\ell_j.
\]
Therefore \eqref{eq:per-error} is $\|H_d-G\|_{L_1}$ and
$I^d=2^{-d}$ is the initial error.

Delete zero coefficients and zero-mass product kernels as in
Remark~\ref{rem:zero-terms}.  Lemma~\ref{lem:product}, with the affinity
bounds in Lemmas~\ref{lem:star-affinity}, \ref{lem:ext-affinity}, and
\ref{lem:per-affinity}, now shows that every individual rule satisfying
\[
 E^{\mathcal C}_{N,d}(P,a)
 \le\varepsilon E^{\mathcal C}_{0,d}
\]
must satisfy
\begin{equation}\label{eq:individual-bound}
 N\ge\frac{(1-\varepsilon)^2}{1+\varepsilon}
 q_{\mathcal C}^d,
\end{equation}
where
\begin{equation}\label{eq:q-values}
 q_{\mathrm{star}}=\frac12+\frac1{\sqrt3},\qquad
 q_{\mathrm{ext}}=\frac{75(9+4\sqrt2)}{784},\qquad
 q_{\mathrm{per}}=\frac{25}{16}.
\end{equation}

We finally pass from individual rules to the infimum in
\eqref{eq:complexity}.  Suppose
\[
 \disc^{\mathcal C,+}_1(N,d)
 \le\varepsilon E^{\mathcal C}_{0,d}.
\]
For every $\eta>0$ with $\varepsilon+\eta<1$, the definition of the
infimum gives a rule with error at most
\((\varepsilon+\eta)E^{\mathcal C}_{0,d}\).  Applying
\eqref{eq:individual-bound} with $\varepsilon+\eta$ and letting
$\eta\downarrow0$ yields
\[
 N\ge\frac{(1-\varepsilon)^2}{1+\varepsilon}q_{\mathcal C}^d.
\]
This proves \eqref{eq:main-star}--\eqref{eq:main-per}.  The QMC estimates
follow from \eqref{eq:qmc-subclass}.

Each number in \eqref{eq:q-values} is larger than one.  Therefore, for every
fixed $0<\varepsilon<1$ and every $1<c<q_{\mathcal C}$, the constant factor
in \eqref{eq:individual-bound} is eventually dominated by
\((q_{\mathcal C}/c)^d\).  Hence
\(N^{\mathcal C,+}_1(\varepsilon,d)\ge c^d\)
for all sufficiently large $d$.  This proves the final assertion of
Theorem~\ref{thm:main}.

\begin{remark}[Why nonnegative weights are essential]
The proof uses
\[
 \sqrt{\sum_i b_iF_i}\le\sum_i\sqrt{b_iF_i},
 \qquad b_iF_i\ge0.
\]
There is no corresponding estimate for a signed sum that retains the mass
control used in \eqref{eq:B-upper}. Thus the paper makes no claim for
arbitrary real coefficients.
\end{remark}

\section*{Acknowledgement}

The work of J.~D. was supported in part by a grant of access to OpenAI
models through the ChatGPT for Academic Researchers program.

\section*{Declaration of generative AI use}

J.~D. used ChatGPT during exploratory development, literature search,
and the preparation of the exposition and LaTeX source.
The authors take full responsibility for all mathematical arguments,
calculations, references, and conclusions in this paper.

\end{document}